\documentclass[]{interact}

\usepackage{enumerate}
\usepackage{epstopdf}
\usepackage[caption=false]{subfig}
\usepackage[numbers,sort&compress]{natbib}

\makeatletter
\def\NAT@def@citea{\def\@citea{\NAT@separator}}
\makeatother

\makeatletter
\renewcommand*\env@matrix[1][*\c@MaxMatrixCols c]{%
  \hskip -\arraycolsep
  \let\@ifnextchar\new@ifnextchar
  \array{#1}}
\makeatother

\theoremstyle{plain} 
\newtheorem{theorem}{Theorem}[section]

\newtheorem{conclusion}[theorem]{Conclusion}

\theoremstyle{definition} 

\newtheorem{remark}[theorem]{Remark}

\newtheorem*{assumption*}{Assumption}

\usepackage{comment}
\usepackage{arydshln}
\usepackage{soul}
\usepackage[color=yellow!30, linecolor=orange, colorinlistoftodos]{todonotes}
\usepackage[breakable]{tcolorbox}
\tcbuselibrary{theorems}

\colorlet{mycyan}{cyan!20}
\colorlet{myorange}{orange!40}

\colorlet{myyellow}{yellow!40}
\sethlcolor{myyellow}

\usepackage{amsmath}
\usepackage{amssymb}
\usepackage{abraces}
\usepackage[utopia]{mathdesign}
\usepackage{hyperref}
\usepackage{tikz}
\usetikzlibrary{matrix, calc}
\usetikzlibrary{positioning}
\usepackage{datetime}

\allowdisplaybreaks

\makeatletter
\def\widebreve{\mathpalette\wide@breve}
\def\wide@breve#1#2{\sbox\z@{$#1#2$}%
     \mathop{\vbox{\m@th\ialign{##\crcr
\kern0.08em\brevefill#1{0.8\wd\z@}\crcr\noalign{\nointerlineskip}%
                    $\hss#1#2\hss$\crcr}}}\limits}
\def\brevefill#1#2{$\m@th\sbox\tw@{$#1($}%
  \hss\resizebox{#2}{\wd\tw@}{\rotatebox[origin=c]{90}{\upshape(}}\hss$}

\begin{document}

\title{Discretization of the wave equation on a metric graph}

\author{
\name{Sergei A. Avdonin\textsuperscript{1}
\thanks{
CONTACT Sergei~A.~Avdonin.  Email: saavdonin@alaska.edu}, 
Aleksander S. Mikhaylov\textsuperscript{2,3},
Victor S. Mikhaylov\textsuperscript{2}.
Abdon E. Choque-Rivero\textsuperscript{4,}}
\affil{\textsuperscript{1} 
Department of Mathematics \& Statistics, 
University of Alaska Fairbanks, Fairbanks,
Alaska, USA;
\textsuperscript{2}~St.Petersburg Department of Steklov Mathematical Institute, 
Russian Academy of Sciences, St. Petersburg, Russia;
\textsuperscript{3}~Department of Mathematics, 
St. Petersburg State University, St. Petersburg, Russia;
\textsuperscript{4}~Institute of Physics and Mathematics, 
Universidad Michoacana de San Nicolas de Hidalgo, Michoacan, Mexico}
}

\maketitle

\begin{abstract}
   The question of what conditions should be set at the nodes of a discrete graph for the wave equation with discrete time is investigated. The variational method for the derivation of these conditions is used. A parallel with the continuous case is also drawn. As an example the problem of shape controllability from the boundary is studied. 
\end{abstract}

\begin{keywords}
Discrete graph, node conditions, wave equation.
\end{keywords}
 
\begin{amscode}
35L05, 35Q93, 93B05, 93C20, 37M15, 37J51
 \end{amscode}
\section{Introduction}\label{sec1}
Under differential equation networks (DENs) or, in other words,  quantum graphs we understand differential equations on metric graphs coupled by certain vertex matching conditions. These models play a fundamental role in many problems of
	science and engineering. The range for applications  of DENs is enormous and continues to grow;  we will mention a few of them.
	
	-- {\it Structural Health Monitoring.} DENs, classically, arise in the study of stability, health, and oscillations of 
	flexible structures that are made of strings, beams, cables, and struts
 \cite{LLS-book, ref15, KIIK, ref17}

	-- {\it Water, Electricity, Gas, and Traffic Networks.} An important example of DENs is the Saint-Venant system
	of equations, which model hydraulic networks for water supply and irrigation 
  \cite{GL-journal}
 and first-order hyperbolic equations \cite{ref1,ref2,ref3,ref4,ref5,ref6} and the
	isothermal Euler equations for describing the gas flow
	through pipelines \cite{BCd-book, Hante-book}.
	Other important
	examples of DENs include the telegrapher equation for modeling electric networks \cite{ABG-journal}, 
	the diffusion equations in power networks \cite{ref14}, and
	Aw-Rascle equations for describing road traffic dynamics \cite{CGHS-journal}, see also \cite{ref7}
 for traffic flow on networks and \cite{ref9,ref10,ref11} 
 for modeling groundwater flow.
	
	-- {\it Nanoelectronics and Quantum Computing.} Mesoscopic quasi-one-dimensional structures  such as
	quantum, atomic, and molecular wires are the subject of extensive experimental and theoretical studies \cite{H-book, JR-book}
 The simplest model describing conduction in quantum wires is the Schr\"odinger operator
	on a planar graph. For similar models appear in
	nanoelectronics, high-temperature superconductors, quantum computing and studies of quantum chaos see 
 \cite{MP-journal, ref12, KS1-journal}.

	-- {\it Material Science.} Quantum graphs arise in analyzing hierarchical materials like ceramic and metallic foams, percolation
	networks, carbon and graphene nano-tubes, and graphene ribbons  
 \cite{Hirt1974, ref8}.

	-- {\it Biology.} Challenging problems involving ordinary and partial
	differential equations on graphs arise in 
	signal propagation in dendritic trees, particle dispersal in respiratory systems, species
	persistence and biochemical diffusion in delta river systems 
 \cite{BC-journal, AB-journal}.

	-- {\it Social Networks.} Examples of social applications are modeling of international trades and space-temporal patterns of information spread
  \cite{ref13}.
	
	There are many papers in the literature devoted to study of the spectral properties of differential operators on graphs and well-posedness of the initial boundary value problems for differential equations on graphs
	and regularity of their solutions \cite{ref19,ref20,ref21,ref22,ref23,ref24,BK-book}.
 On the other hand, numerical methods for solving ODEs and PDEs on graphs have been mostly developed only for very specific problems \cite{ref6,ref25,ref26,ref27,ref28,ref29,ref30,ref31}.
 Only recently more general investigations appeared directed to developing the finite element method for elliptic and parabolic equations 
    on graphs \cite{ref32, ref33, ref34}.
	
	In the present paper we study the problems of discretization of the wave equation on metric graphs. The main attention is put to discretization of the Kirchoff -- Neumann matching conditions at the internal vertices. This important problems has not got a proper attention in the literature. To derive the matching conditions we apply Hamilton's principle and variational methods in both continuous and discrete cases. Then we discuss a proper choice of the nodal weights for the discrete model to insure the same reflection and transmission coefficients as the original continuous model has. The paper is organized as follows. In Section 2 we derive the matching conditions for metric graphs and in Section 3 --- for discrete graphs. In Section 4 we discuss transmission and reflection of the wave at the internal nodes in both, continuous and discrete, situations. In Section 5 we solve direct and control problems for the wave equation on a discrete star graph.

\section{Conditions at nodes on a metric graph}\label{sec2}
Let $ \Omega $ be a finite connected compact graph. The graph consists of edges $E = \{e_1, \ldots, e_{N}\}$ connected at vertices $V = \{v_1 \ldots, v_{M}\}$. Each edge $ e_j \in E $ is identified with an interval $ (0, l_j) $ of the real line. The boundary $ \Gamma = \{v_1, \ldots, v_m \} $ of $ \Omega $ is the set of vertices whose degree is equal to one (outer nodes). Let $E(v)$ be a set  of edges incident to $v$.
	In what follows, we assume that some of the boundary vertices are clamped, and the non-homogeneous Dirichlet boundary condition is imposed on the other part. 
	
	
	The space of real square-integrable functions on the graph $ \Omega $ is denoted by $L_2(\Omega ):=
	\bigoplus_{i=1}^{N}L_2(e_i).$ For the function $u \in L_2(\Omega)$ we will write
	\begin{equation*}
		u:= \left\{ u^i\right\}_{i=1}^N\ ,\ u^i\in L_2(e_i).
	\end{equation*}
	
\begin{figure}[h]
		\includegraphics[width=0.9\textwidth]{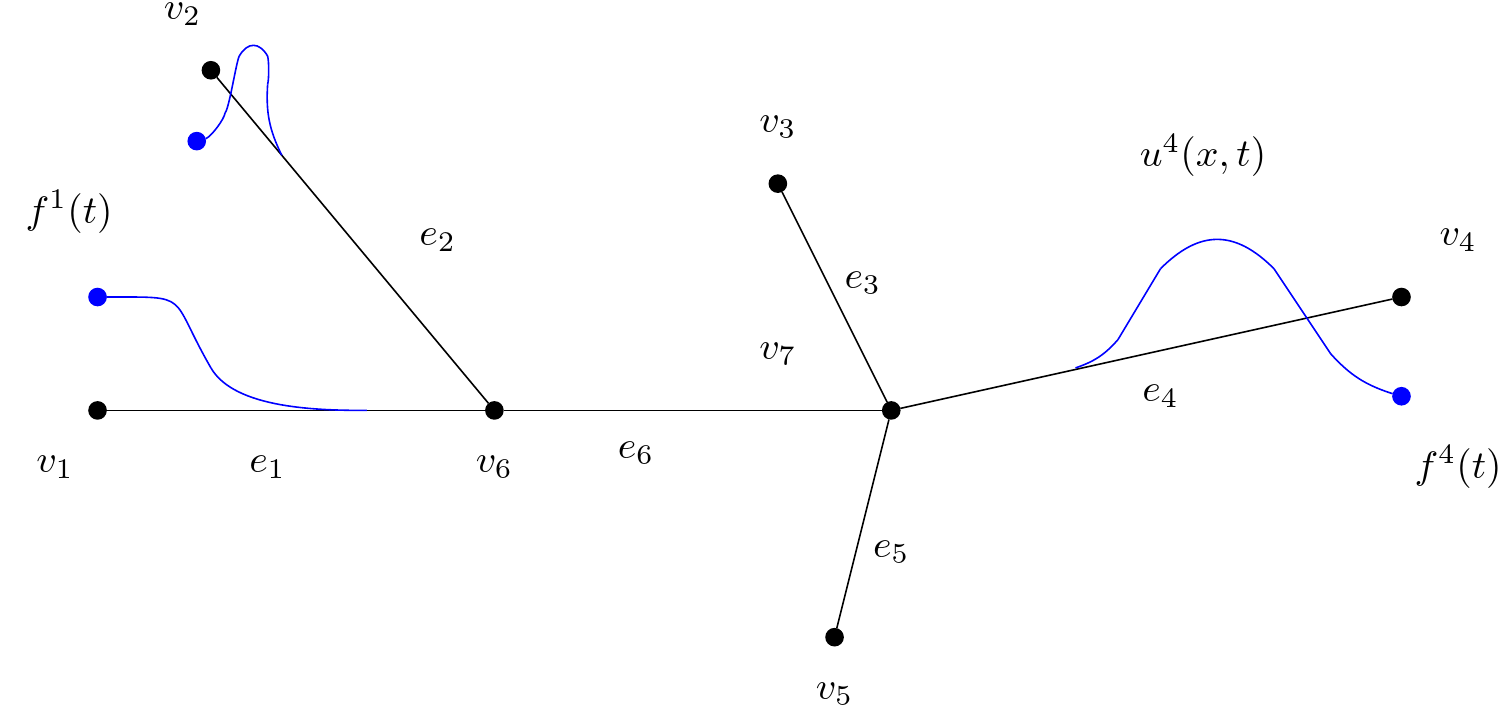}
		\caption{Wave equation on a metric graph with six edges and seven vertices.
  On the edge $e_1$ (resp. $e_4$) the function $f^1(t)$ (resp. $f^4(t)$) is acting.
   On the edge $e_4$ the corresponding displacement $u(x,t)$ is represented.
  \label{figure1}}		
	\end{figure}

	By $ \textbf{C} $ we denote the direct problem on the metric  graph:
	\begin{align}
		u^i_{tt}(x,t)-u^i_{xx}(x,t)=0 \quad x \in e_i, \label{Wave_eqn} \\
		u^i(v,t)=u^j(v,t),  \quad e_i,\ e_j\in E(v),\quad v\in V\backslash\Gamma,  \label{Cont}\\ 
		\sum_{ e_i\in E(v)}\partial u^i(v,t)=0,         v\in V\backslash \Gamma, \label{Kirch}\\
		u(v,t) = f(t)\quad \text{for}\,\, v\in\Gamma, \label{Bound}\\
		u(x,0) = 0, u_t(x,0) = 0, \ \quad \text{for}\,\, x\in\Omega. \label{Init}
	\end{align}
	Here $\partial u^i(v,t)$ is the derivative of $u^i$ at $v$ along $e_i$ in the direction away from $v$.
	This system describes small transverse oscillations of the graph (Figure  \ref{figure1}), and $u(x,t)$ denotes the displacement of the point $x\in\Omega$ of the graph at time $t$ from the equilibrium position. It is assumed that the boundary of the graph moves according to the law $f(t)$, and in the initial moment of time  the graph was at rest in an equilibrium state.

	In order to write a similar system of equations  for a discrete graph, let us recall how equations (\ref{Wave_eqn})--(\ref{Init}) arose. We introduce the kinetic and potential energies:
	\begin{align}
		T(t)=\frac{1}{2}\int_\Omega u_t^2(x,t)\,dx, \label{Kin_en} \\
		U(t)=\frac{1}{2}\int_\Omega u_x^2(x,t)\,dx. \label{Pot_en} 
	\end{align}
	According to the principles of classical mechanics (Hamilton's principle of least action), the system passes from state 1 at time $t_1$ to state 2 at time $t_2$ in such a way that the variation of the action functional vanishes
	\begin{equation}
		S[u]= \int_{t_1}^{t_2} L(t)\,dt,\quad L=T-U\text{ (Lagrangian)},
	\end{equation}
	along the true trajectory of movement. It means that  
	$$
	\delta S[u,h]=0 \quad \text{for all}\quad h: h|_{t=t_1}=0, h|_{t=t_2}=0, h|_\Gamma=0.
	$$
	Calculating the variation of the functional and integrating by parts, we get the following:
	\begin{align*}
		0=&\delta S[u,h]=\int_{t_1}^{t_2}\int_{\Omega} (u_t h_t -u_x h_x)\,dx\,dt\\
		 =&-\int_{t_1}^{t_2}\int_{\Omega} (u_{tt} -u_{xx})h\,dx\,dt+ \int_{t_1}^{t_2}\int_{V\backslash\Gamma} \partial u h \,dx\,dt,
	\end{align*}
	where the integrals over the boundary $\Gamma$ and for $t=t_{1,2}$ vanish because of the requirements on $h$. From the last equality and the main Lemma of the Calculus of Variations (the Dubois-Reymond Lemma), we obtain the equation (\ref{Wave_eqn}) and the condition at the internal node (\ref{Kirch}), the so-called Kirchhoff condition. The condition (\ref{Cont}) is a natural continuity condition. The first condition in (\ref{Init}) means that we know the position of the system at time $t=t_1=0$, and instead of writing that we know the position of the system at time $t=t_2$ we write the second condition in (\ref{Init}). Thus, the second initial condition does not follow from the principle of least action, but replaces the condition of knowing the position of the system at the final moment of time.
	\begin{remark}
		If there is a point mass $m_j$ at the node $v_j\in V$, $j=m+1,\ldots,N+1$, then the kinetic energy of this mass $m_j {\displaystyle\frac{u^2_t(v_j,t)}{2}}$ must be added to the kinetic energy $T$. Then the Lagrangian will change and the variation of the action functional will look like this:
		$$
		\delta S[u,h]=\ldots+\int_{t_1}^{t_2}m_j u_t(v_j,t) h_t\,dt\ =\ldots -\int_{t_1}^{t_2}m_j u_{tt}(v_j,t) h\,dt.
		$$
		The last term contributes to the knot, and the Kirchhoff condition takes the following form
		$$
		m_ju_{tt}(v_j,t)=\sum_{e_i\in E( v_j)}\partial u^i(v_j,t).
		$$ 
	\end{remark}
	This equality is nothing but Newton's second law - the product of mass and acceleration is the sum of all forces acting on the body. On the right side, just $\partial u^i(v_j,t)$ is the tension force from the side of the edge $e^i$ acting on the node $v_j$.
	
	In the next section we consider discrete analogue of system (\ref{Wave_eqn})--(\ref{Init})  -- a discrete wave equation on a discrete graph.

\section{Conditions at nodes on a discrete graph}
	
\begin{figure}[h]
		\includegraphics[width=0.9\textwidth]{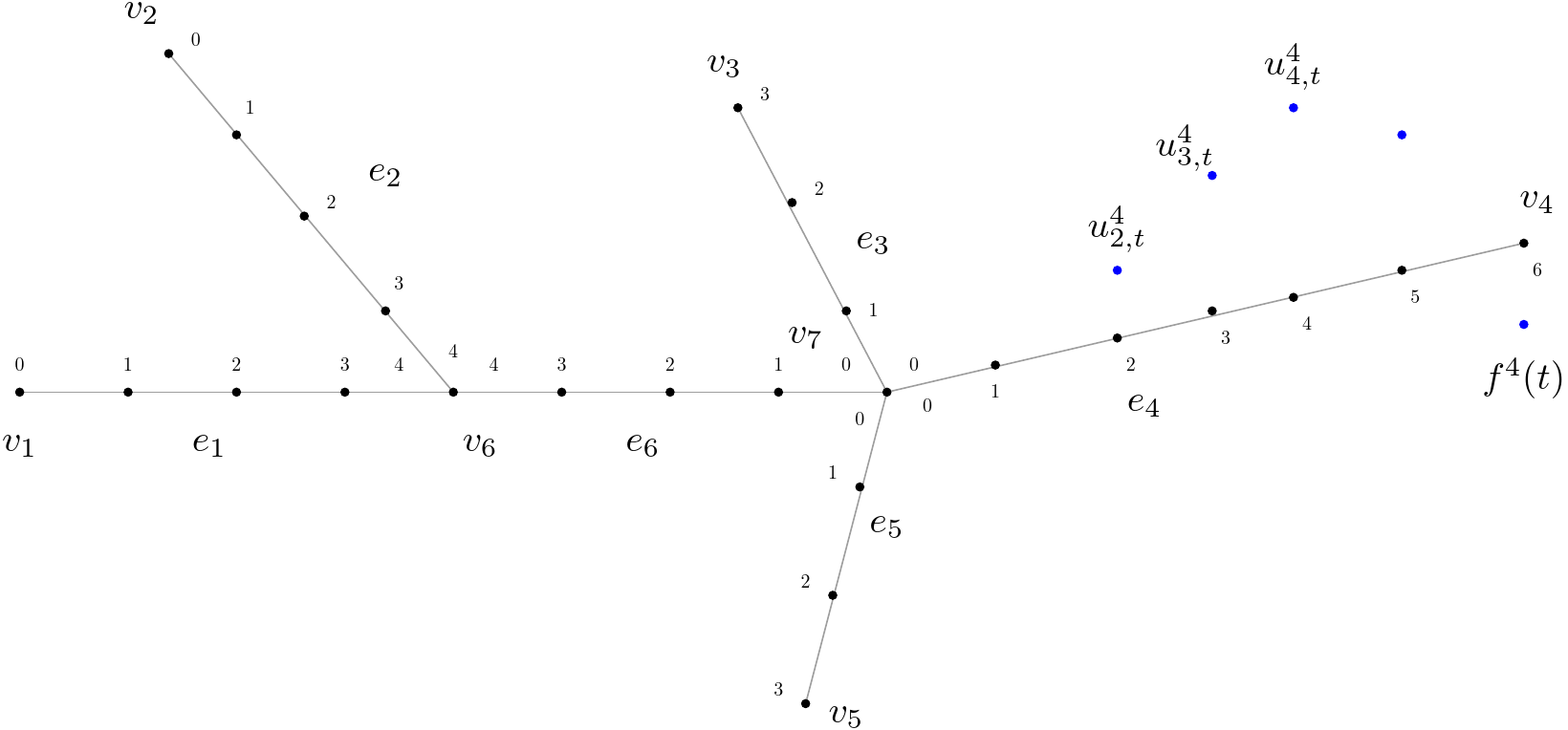}
		\caption{Discrete Wave equation on a discrete graph with six edges and seven vertices. On the edge $e_4$ is acting the function $f^4(t)$. On the same edge
  the discrete displacement of $u^4_{j,t}$ is represented.\label{figure2}}
	\end{figure}

	In this section, we assume that our metric graph $\Omega$ has been replaced by a discrete graph $\Omega_D$. This  means that the vertices have remained in their places, and each edge has become a finite discrete set of points, or, in other words, each edge $e_j\in E$ is now identified not with the interval $(0,l_j)$ of the real line, but with a finite set of numbers $a^j_0<a^j_1<\ldots<a^j_{N_j}$. For simplicity, we can assume that $a^j_0=0$, $a^j_1=1$,\ldots$a^j_{N_j}=N_j$ (Figure \ref{figure2}). 

	In addition, we assume that the time $t$ is also discrete, i.e. we replace the interval $(0,T)$ with the discrete set $t_0<t_1<\dots< t_T$ or, for simplicity, $0<1<\ldots <T$.
	The space of real square summable functions on the graph $ \Omega_D $ is denoted by $L(\Omega_D):=
	\bigoplus_{i=1}^{N}\mathbb{R}^{N_i}.$  For the function $u \in L(\Omega_D)$ we will write
	\begin{equation*}
		u:= \left\{ u^i\right\}_{i=1}^N,\quad u^i=(u^i_0,\ldots,u^i_{N_i}),\quad u^i_j\in \mathbb{R},\ j=0,\ldots,N_i.
	\end{equation*}
	
	In order to unify the value at the node $v_k$, we introduce the function $F(k)=\left[ \begin{array}l 0\\ N_k\end{array}\right.$ such that
	$$
	u(v_k)=u^{j(k)}_{F(k)},\quad{\text{where } j(k)\in E(v_k)},\quad k=1,\ldots,M.
	$$
	
	Since we want to consider a discrete wave equation, the  function $u$ must depend on two variables - coordinate and time. Time is now discrete, so we introduce this dependence through an additional index on the vector $u^i$. We will write 
	\begin{equation*}
		u(t):= \left\{ u^i(t)\right\}_{i=1}^N,\quad u^i(t)=(u^i_{0,t},\ldots,u^i_{{N_i},t}),\quad u^i_{j,t}\in \mathbb{R},\ j=0,\ldots,N_i.
	\end{equation*}
	Before reformulating $ \textbf{C} $ -- a direct problem on a metric graph and writing a similar discrete one, we consider the discrete analogs of the kinetic and potential energy functionals:
	\begin{align}
		T_D(t)=&\sum_{i=1}^N T^e_i(t)+\sum_{i=1}^{M}T^v_i(t)=
		\sum_{i=1}^N \sum_{j=1}^{N_i-1}\frac{(u^i_{j,t}-u^i_{j,t-1} )^2}{2}  +\sum_{i=1}^{M} \frac{(u^{j(i)}_{F(i),t}-u^{j(i)}_{F(i),t-1} )^2}{2} \label{Kin_en_d},
\\
		U_D(t)=&\sum_{i=1}^N U^e_i(t)= \sum_{i=1}^N \sum_{j=1}^{N_i} \frac{(u^i_{j,t}-u^i_{j-1,t} )^2}{2}. \label{Pot_en_d} 
	\end{align}
	Here $T^e_i$, $U^e_i$, $i=1,\ldots N$ are the kinetic and potential energies from the interior points of the $e_i$ edges, and $T^v_i$, $i=1,\ldots N+1$ --- kinetic energy from $v_i$ vertices. In order to write the action functional in the discrete  case, it was necessary to integrate the Lagrangian over the interval $(t_1,t_2)$. In a discrete situation, the formula looks like this:
	$$
	S[u]=\sum_{t=1}^{T}T_D(t)- \sum_{t=0}^{T}U_D(t).
	$$
	
	Let us apply the principle of least action  to this quadratic functional. According to this principal the variation of the action functional $S[u]$ along the true trajectory of motion vanishes:
	
	\begin{align*}
		0=&\delta S[u,h]= \sum_{t=1}^{T}\sum_{i=1}^N \sum_{j=1}^{N_i-1}(u^i_{j,t}-u^i_{j,t-1} ) (h^i_{j,t}-h^i_{j,t-1} )+\sum_{t=1}^{T}\sum_{i=1}^{M}(u^{j(i)}_{F(i),t}-u^{j(i)}_{F(i),t-1} )(h^{j(i)}_{F(i),t}-h^{j(i)}_{F(i),t-1} )-
		\\
		&- \sum_{t=0}^{T}\sum_{i=1}^N \sum_{j=1}^{N_i} (u^i_{j,t}-u^i_{j-1,t} )(h^i_{j,t}-h^i_{j-1,t} )
		\\
		=&\sum_{t=1}^{T}\sum_{i=1}^N \sum_{j=1}^{N_i-1} h^i_{j,t}(u^i_{j,t}-u^i_{j,t-1} - (u^i_{j,t+1}-u^i_{j,t}))+
		\sum_{t=1}^{T}\sum_{i=1}^{M} h^{j(i)}_{F(i),t}(u^{j(i)}_{F(i),t}-u^{j(i)}_{F(i),t-1} - (u^{j(i)}_{F(i),t+1}-u^{j(i)}_{F(i),t}))\\
		&-\sum_{t=0}^{T}\sum_{i=1}^N \sum_{j=1}^{N_i-1}h^i_{j,t}(u^i_{j,t}-u^i_{j-1,t} -(u^i_{j+1,t}-u^i_{j,t}))-
		\sum_{t=0}^{T}\sum_{i=1}^N h^i_{0,t}(u^i_{0,t}-u^i_{1,t}) -\sum_{t=0}^{T}\sum_{i=1}^N h^i_{N_i,t}(u^i_{N_i,t}-u^i_{N_i-1,t}).
	\end{align*}
	
	Now we use the fact that the factors $h^i_{j,t}$ in the above sum are  arbitrary except boundary $\Gamma$  and, in view of this, we obtain the equation at the internal points of the edges:
	\begin{align}
		0=&(u^i_{j,t}-u^i_{j,t-1} - (u^i_{j,t+1}-u^i_{j,t}))-(u^i_{j,t}-u^i_{j-1,t} -(u^i_{j+1,t}-u^i_{j,t}))\nonumber\\
		=&-(u^i_{j,t+1}+u^i_{j,t-1}-u^i_{j+1,t}-u^i_{j-1,t}), \label{Wave_eqn_d}
	\end{align}
	where $i=1,\ldots,N$, $t=1,\ldots,T-1$, $j=1,\ldots,N_i-1$.
	At the boundary points we assumed the same conditions as in (\ref{Bound}):
	$$
	u(v_k,t) = f^k(t)\quad \text{for}\,\, v_k\in\Gamma. 
	$$
	It remains to write the condition at the internal nodes. Consider vertex $v_k$. 
	Collecting together the terms with the factor $h^{j(k)}_{F(k),t}$ and equating it to zero, we get:
	\begin{equation}
		(u^{j(k)}_{F(k),t}-u^{j(k)}_{F(k),t-1} - (u^{j(k)}_{F(k),t+1}-u^{j(k)}_{F(k),t})) - \sum_{i: e_i\in E(v_k)} (u^i_{F(i),t}-u^i_{|F(i)-1|,t})=0\label{Kirch_d},
	\end{equation}
	
	where  $|F(i)-1|=\left[ \begin{array}l 1\\ N_i-1\end{array}\right.$ and thus $(u^i_{F(i),t}-u^i_{|F(i)-1|,t})$ is a discrete analogue of operation $\partial$ on a metric graph.

	Assume that the number of elements in $E(v_k)$ is $p_k$.
	Then the last equality is transformed to
	\begin{equation}
		-u^k_{F(k),t+1}-u^k_{F(k),t-1} +(2-p_k)u^k_{F(k),t} + \sum_{e_i\in E(k)}u^{i}_{|F(i)-1|,t}=0.\label{Kirch_d_1}
	\end{equation}
	Thus, we have obtained a discrete wave equation on the graph (\ref{Wave_eqn_d}) and conditions at the nodes (\ref{Kirch_d}) or (\ref{Kirch_d_1}), which are analogues of (\ref{Wave_eqn}) and (\ref{Kirch})
	\begin{remark}
		If at the node $v_k\in V$, there is a point mass $m_k$, then to the kinetic energy $T_D$ it is necessary to add the kinetic energy of this mass $m_k {\displaystyle\frac{(u(v_k,t)-u(v_k,t-1 ))^2}{2}}$. Then the Lagrangian will change and the variation of the action functional and the condition at the node will look like this:
		\begin{equation}
			-(1+m_k)u^k_{F(k),t+1}-(1+m_k)u^k_{F(k),t-1} +(2+2m_k-p_k)u^k_{F(k),t} + \sum_{e_i\in E(v_k)}u^{i}_{|F(i)-1|,t}=0.\label{Kirch_d_1m}
		\end{equation}
	\end{remark}
	\begin{remark}
		If we put $m_k=-1$ then (\ref{Kirch_d_1m}) becomes
		\begin{equation}
			\sum_{e_i\in E(v_k)}(u^{i}_{|F(i)-1|,t}-u^k_{F(k),t})=0,\label{Kirch_d_1m0}
		\end{equation}
		which is a direct discrete analog of the Kirchhoff condition (\ref{Kirch}). Note that the formal choice of $m_k=-1$ practically means that there is zero mass at node $v_k$, because by our definition, the kinetic energy functional (\ref{Kin_en_d}) contained the energy $T^v_k$ from the unit mass at the node $v_k$.
	\end{remark}
	\begin{remark}
		If we put \begin{equation}
			m_k=\frac{p_k-2}{2},\label{mass_cond}
		\end{equation}
		then (\ref{Kirch_d_1m}) become
		\begin{equation}
			-\frac{p_k}{2}u^k_{F(k),t+1}-\frac{p_k}{2}u^k_{F(k),t-1}  + \sum_{e_i\in E(v_k)}u^{i}_{|F(i)-1|,t}=0.\label{Kirch_d_1mm}
		\end{equation}	
	\end{remark}
	We are faced with the choice of which of the conditions at the nodes - (\ref{Kirch_d_1}), (\ref{Kirch_d_1m0}) or (\ref{Kirch_d_1mm}) to prefer in the mathematical study of the direct and inverse problems. From a mathematical point of view, we can take any of these conditions, because in the limit, when a continuous equation is approximated by a discrete one, they will turn into the Kirchhoff conditions. It seems that (\ref{Kirch_d_1m0}) looks more preferable than the others, because it is already a direct analogue of the Kirchhoff condition (\ref{Kirch}). In the next section, we will study this issue and come to the conclusion that the correct condition is actually the condition $(\ref{Kirch_d_1mm})$. For example we will demonstrate the following fact
	\begin{remark}
		In the case when two edges $e_1$ and $e_{i_2}$ converge at the vertex $v_k$ and, consequently, $p_k=2$, from the formula (\ref{mass_cond}) we get that $m_k=0$, and ( \ref{Kirch_d_1mm}) will become
		$$
		-u^1_{F(1),t+1}-u^1_{F(1),t-1}  + u^1_{|F(1)-1|,t} + u^{i_2}_{|F(i_2)-1|,t}=0.
		$$
		The last equation is the same as (\ref{Wave_eqn_d}), which in turn means that these two edges $e_1$ and $e_{i_2}$
		can be replaced by one by removing the $v_k$ vertex from the graph.
		
		We can conclude that the choice of the mass $m_k$ according to the formula (\ref{mass_cond}) is more correct than, for example, the option when we choose $m_k=-1$.
	\end{remark}

\section{Transmission and reflection of waves at the internal node}
	
	In this section, as a model situation, we will consider a star graph\, --\, a connected graph in which all edges are incident to one  internal vertex. A star-graph with $k+1$ vertices is usually denoted $S_k$, where $k$ is called the order of the star. In our notation, all edges $e_1,\ldots,e_k$ are incident to the vertex $v_{k+1}$, the center of the star-graph. Each edge $e_i$ starts at $v_i$ $i=1,\ldots,k$ and ends at $v_{k+1}$. 
	
	{\noindent\bf Metric graph.}
	
	To understand  the processes of transmission and reflection we recall what happens in the case of a metric graph for the problem $ \textbf{C} $ when the boundary conditions (\ref{Bound}) have the form:
	$$
	u^1(0,t)=\delta(t),\quad u^i(0,t)=0,\quad i=2,\ldots,k.
	$$
	
	This statement corresponds to the situation when  all but one of the boundary vertices are clamped and the first vertex  is clicked at the initial moment of time. Due to the fact that on each edge we have the wave equation (\ref{Wave_eqn}), this wave will propagate inside the graph along the edge $e_1$ according to the rule $u^1(x,t)=\delta(t-x)$ until it reaches $v_{k+1}$. Since the length of the edge $e_1$ is equal to $l_1$, and the speed of propagation is equal to $1$, then at time $t=l_1$ the delta function reaches  the vertex $v_{k+1}$, the center of the star-graph. At this moment of time, the conditions at the node come into force -- the continuity conditions (\ref{Cont}) and the Kirchhoff condition (\ref{Kirch}).
	Using them, one can see that $\frac{2}{k}\delta$ will run along each outgoing edge $e_i$ $i=2,\ldots,k$, and the reflected wave will run along the edge $e_1$ in the opposite direction $-\frac{k-2}{k}\delta$. Namely:
	$$
	u^1(x,t)=-\frac{k-2}{k}\delta(t+x-2l_1), \quad u^i(x,t)=\frac{2}{k}\delta(t+x-l_1-l_i),\quad i=2,\ldots,k,
	$$
	where $l_1<t<l_1+l$, $l:=\min\{l_1,\ldots,l_k\}$. Indeed, the continuity condition is satisfied, because for $x=l_1$ and $t=l_1$ on the first edge the solution looks like the sum of the incident and reflected waves, and on the other edges as the transmitted one:
	$$
	u^1=\delta(t-x)-\frac{k-2}{k}\delta(t+x-2l_1), \quad u^i=\frac{2}{k}\delta(t+x-l_1-l_i),\quad  i=2,\ldots,k.
	$$
	The continuity and the Kirchhoff conditions (\ref{Cont}), (\ref{Kirch}) are easily follow from representations above.

	It will be important to us  that the delta function falling on the node passes through it and runs further along the outgoing edges with the coefficient $2/k$ and is reflected and runs in the opposite direction along the incoming edge with the coefficient $-(k-2)/k $.
	
	We want to observe exactly the same picture in the situation when the graph and equations have become discrete.
	
	{\noindent\bf Discrete graph.}
	
	Let $k=3$ and all three edges of the star graph have the same length, i.e. are parameterized by the set of numbers $0,1,\ldots,N$. Recall that time is discrete and the dependence on it is taken into account through the second subscript. So  function under consideration has a form:
	
	\begin{equation*}
		u:= \left\{ u^i\right\}_{i=1}^3,\quad u^i=(u^i_{0,t},\ldots,u^i_{N,t}),\quad u^i_{j,t}\in \mathbb{R},\ j=0,\ldots,N.
	\end{equation*}
	For $i=1,2,3$, the index value $j=0$ corresponds to the graph vertex $v_i$, and the index value $j=N$ corresponds to the vertex $v_4$ --- the center of the star. This function satisfies the equation (\ref{Wave_eqn_d}), continuity, initial and boundary conditions:
	\begin{align}
		u^i_{j,t+1}+u^i_{j,t-1}-u^i_{j+1,t}-u^i_{j-1,t}=0,\   j=1,\ldots,N-1, t=0,\ldots\label{Wave_eqn_d3}\\
		u^1_{N,t}=u^2_{N,t}=u^3_{N,t},\  t=1,\ldots \label{Cont_d3}\\
		u^1_{0,0}=1,\ u^2_{0,0}= u^3_{0,0}=u^i_{0,t}=0,\ i=1,2,3,\ t=1,\ldots  \label{Bound_d3}\\
		u^i_{j,0}=u^i_{j,-1}=0,\ i=1,2,3,\ j=1,\ldots,N.\label{Init_d3}
	\end{align}
	In the last line we formally introduced the initial condition for  $t=-1$. This is done to write an analog of the condition $u_t(x,0)=0$.
	
	To this system one must add a condition at the internal node (the analogous to the Kirchhoff condition). As mentioned in the previous section, several variants of this condition are possible. Let us take a look at each of the following three conditions:
	\begin{align}
		(u^1_{N-1,t}-u^1_{N,t})+(u^2_{N-1,t}-u^2_{N,t})+(u^3_{N-1,t}-u^3_{N,t})=0,\label{Kirch_d_1m03}\\
		-u^1_{N,t+1}-u^1_{N,t-1} -u^1_{N,t} + u^1_{N-1,t}+u^2_{N-1,t}+u^3_{N-1,t}=0,\label{Kirch_d_13}\\
		-\frac{3}{2}u^1_{N,t+1}-\frac{3}{2}u^1_{N,t-1}  + u^1_{N-1,t}+u^2_{N-1,t}+u^3_{N-1,t}=0.\label{Kirch_d_1mm3}
	\end{align}
	
	We present the solution of each of the three problems in tables. For simplicity we assume that $N=3$. For the
	problem (\ref{Wave_eqn_d3})--(\ref{Init_d3}) and condition (\ref{Kirch_d_1m03}) we get
	
	\begin{equation*}
		\begin{tabular}{c|ccc|c|ccl}
			6  & \bf0 & -1/3 & 0 & 0 & 0 & -1/3 & \bf0 \\
			5  & \bf0 & -2/3 & 0 & 0 & 0 & 1/3 &\bf 0 \\
			4  & \bf0 & 1/3 & -2/3 & 0 & 1/3 & 1/3 & \bf0 \\
			3  & \bf0 & 0 & 1/3 & 1/3 & 1/3 & 0 & \bf0 \\
			2  & \bf0 & 0 & 1 & 1/3 & 0 & 0 & \bf0 \\
			1  & \bf0 & 1 & 0 & 0 & 0 & 0 & \bf0 \\
			0  & \bf1 & \bf0 & \bf0 & \bf0 & \bf0 & \bf0 & \bf0 \\
			-1 & \bf0 & \bf0 & \bf0 & \bf0 & \bf0 & \bf0 & \bf0 \\
			\hline 
			t/j & 0 & 1 & 2 & 3 & 2 & 1& 0\\
			&i&=&1& &i&=&2,3
		\end{tabular}
	\end{equation*}

	For the condition (\ref{Kirch_d_13}) we get
	\begin{equation*}
		\begin{tabular}{c|rrr|r|rrl}
			6  & \bf0 & -1 & 2 & -4 & 2 & -1 & \bf0 \\
			5  & \bf0 & 0 & -1 & 2 & -1 & 1 &\bf 0 \\
			4  & \bf0 & 0 & 0 & -1 & 1 & 0 & \bf0 \\
			3  & \bf0 & 0 & 0 & 1 & 0 & 0 & \bf0 \\
			2  & \bf0 & 0 & 1 & 0 & 0 & 0 & \bf0 \\
			1  & \bf0 & 1 & 0 & 0 & 0 & 0 & \bf0 \\
			0  & \bf1 & \bf0 & \bf0 & \bf0 & \bf0 & \bf0 & \bf0 \\
			-1 & \bf0 & \bf0 & \bf0 & \bf0 & \bf0 & \bf0 & \bf0 \\
			\hline 
			t/j & 0 & 1 & 2 & 3 & 2 & 1& 0\\
			&i&=&1& &i&=&2,3
		\end{tabular}
	\end{equation*}

	For the condition (\ref{Kirch_d_1mm3}) we get
	
	\begin{equation}
		\begin{tabular}{c|ccc|c|ccl}
			6  & \bf0 & 0 & 0 & 0 & 0 & 0 & \bf0 \\
			5  & \bf0 & -1/3 & 0 & 0 & 0 & 2/3 &\bf 0 \\
			4  & \bf0 & 0 & -1/3 & 0 & 2/3 & 0 & \bf0 \\
			3  & \bf0 & 0 & 0 & 2/3 & 0 & 0 & \bf0 \\
			2  & \bf0 & 0 & 1 & 0 & 0 & 0 & \bf0 \\
			1  & \bf0 & 1 & 0 & 0 & 0 & 0 & \bf0 \\
			0  & \bf1 & \bf0 & \bf0 & \bf0 & \bf0 & \bf0 & \bf0 \\
			-1 & \bf0 & \bf0 & \bf0 & \bf0 & \bf0 & \bf0 & \bf0 \\
			\hline 
			t/j & 0 & 1 & 2 & 3 & 2 & 1& 0\\
			&i&=&1& &i&=&2,3
		\end{tabular}\label{table3}
	\end{equation}

	Here, boundary and initial conditions are highlighted in bold text. The tables were filled from the bottom to the top with the values of the $u^i_{j,t}$ functions calculated according to the equations (\ref{Wave_eqn_d3})--(\ref{Init_d3}) and the specified conditions (\ref{Kirch_d_1m03})--(\ref{Kirch_d_1mm3}) at the internal node.
	
	In a discrete situation, the analog of the delta function is the vector $(1,0,0,\ldots)$. This is exactly what we see  at the first vertex $v_1$ (in the tables, this is the first column). Note that only in the last case we have got a complete analogue of the continuous situation, namely, that  the delta function falling on the node passes through it and runs along the outgoing edges with a coefficient of $2/3$ and is reflected and runs in the opposite direction along the incoming edge with coefficient $-1/3$. In the first case this statement is true only  ``in the mean", i.e. if we sum the coefficients over two time intervals $t$ and $t+1$. In other words, the support of the initial perturbation ``spreads" and the reflection from the node occurs in ``two cycles". While when choosing the correct third condition, the reflection from the node occurs in ``one cycle" and the support does not ``spread".
	
	In the second case, we see that the resonance begins at the node and such condition does not model continuous situation well.
	
	Summing up, we can say that
	\begin{conclusion}
		The appropriate correct formulation of the direct discrete problem for a discrete graph would be the equations (\ref{Wave_eqn_d3})-(\ref{Init_d3}) and the condition at the node (\ref{Kirch_d_1mm}).
		\begin{eqnarray}
			u^i_{j,t+1}+u^i_{j,t-1}-u^i_{j+1,t}-u^i_{j-1,t}=0,\ 1\leqslant i \leqslant N,\   1\leqslant j \leqslant N_i-1, t\geqslant 0,\label{Wave_eqn_d33}\\
			u^i_{F(i),t}=u^j_{F(j),t},\ e_i,e_j\in E(v_k),\quad   t\geqslant 1,\label{Cont_d33}\\
			u(v_k,t) = f^k(t)\quad \text{for}\,\, v_k\in\Gamma,\quad t\geqslant 1,  \label{Bound_d33}\\
			u^i_{j,0}=0, \quad u^i_{j,-1}=0,\ 1\leqslant i \leqslant N,\   1\leqslant j \leqslant N_i-1,\label{Init_d33}\\
			-\frac{p_k}{2}u^k_{F(k),t+1}-\frac{p_k}{2}u^k_{F(k),t-1}  + \sum_{e_i\in E(v_k)}u^{i}_{|F(i)-1|,t}=0.\label{Kirch_d_1mm33}
		\end{eqnarray}
  In (\ref{Bound_d33}) and (\ref{Kirch_d_1mm33}),  $k=m+1,\ldots,M$.
	\end{conclusion}

	\section{Control Problems of the Wave equations on Discrete Graph}
	
	
	In this section we present  algorithms  solving the direct and control problems  for the discrete wave equation on simple graphs. We start with  a finite interval (finite number of points). Then we construct a solution to the direct problem on a star-shaped graph, and finally we solve a control problem.

	\subsection{Finite interval} Let $\Omega_D$ be a   graph consisting of two vertices $v_1$ and $v_2$ and one
	edge $e$. For discrete problem we identify it with a   finite set of numbers $0,\ldots,N$, the vertex $v_1$ is identified with $x = 0$	and the vertex $v_2$ is identified with $x = N$.
	The control is applied at $x = 0$ and the system takes the form
	\begin{align}
		u_{j,t+1}+u_{j,t-1}-u_{j+1,t}-u_{j-1,t}&=0,\   j=1,\ldots,N-1, t=0,\ldots\label{Wave_eqn_d31}  \\
		u_{0,t} = f_t,\quad u_{N,t} &= 0,\quad \text{for $t=0,\ldots$}\label{B_cond_d31}  \\
		u_{j,0} = 0,\quad u_{j,-1} &= 0, \ \quad \text{for $l=1,\ldots$}\label{I_cond_d31} 
	\end{align}  
	If $T<N$, the solution $u^{f-}$ to (\ref{Wave_eqn_d31})--(\ref{I_cond_d31}) has a form
	$$
	u^{f-}_{j,t}=\left\{
	\begin{array}l0, \quad 0< t < j,\\
		f_{t-j},\quad j\leqslant t.
	\end{array}
	\right.
	$$
 If $T=N$, a reflected wave occurs.
	If $T>N$, the solution  $u^{f-}$ to (\ref{Wave_eqn_d31})--(\ref{I_cond_d31}) has a form
	\begin{equation}
		u^{f-}_{j,t}=\sum_{n=0}^{\left \lfloor{\frac{t-j}{2N}}\right \rfloor} f_{t-j-2nN} - \sum_{n=1}^{\left \lfloor{\frac{t+j}{2N}}\right \rfloor}  f_{t+j-2nN}, \label{solution_left}
	\end{equation}
	where $\lfloor{\cdot}\rfloor$ denotes the floor function.
 In case that the upper index of the second sum of (\ref{solution_left}) is zero, then we set it
 equal to zero.

	Now we consider the case when the control is applied at $x = N$, then  the system takes the form    
	\begin{align}
		u_{j,t+1}+u_{j,t-1}-u_{j+1,t}-u_{j-1,t}&=0,\   j=1,\ldots,N-1, t=0,\ldots\label{Wave_eqn_d32}  \\
		u_{0,t} = 0,\quad u_{N,t} &= g_t,\quad \text{for $t=0,\ldots$}\label{B_cond_d32}  \\
		u_{j,0} = 0,\quad u_{j,-1} &= 0, \ \quad \text{for $l=1,\ldots$}\label{I_cond_d32} 
	\end{align}
	If $T<N$, solution $u^{g+}$ to (\ref{Wave_eqn_d32})--(\ref{I_cond_d32}) has a form
	$$
	u^{g+}_{j,t}=\left\{
	\begin{array}l0, \quad 0< t < j,\\
		g_{t-N+j},\quad j\leqslant t.
	\end{array}
	\right.
	$$
 As in the previous case, if $T=N$, a reflected wave occurs.
	If  $T>N$, the solution  $u^{g+}$ to (\ref{Wave_eqn_d32})--(\ref{I_cond_d32}) has a form
	\begin{equation}
		u^{g+}_{j,t}=\sum_{n=0}^{\left \lfloor{\frac{t+j-N}{2N}}\right \rfloor} g_{t+j-(2n+1)N} - \sum_{n=0}^{\left \lfloor{\frac{t-j-N}{2N}}\right \rfloor}  g_{t-j-(2n+1)N}. \label{solution_right}
	\end{equation}
	
	\subsection{Star-shaped graph}
	Let us consider a three star graph:
	
	\begin{equation*}
		u:= \left\{ u^i\right\}_{i=1}^3,\quad u^i=(u^i_{0,t},\ldots,u^i_{N_i,t}),\quad u^i_{j,t}\in \mathbb{R},\ j=0,\ldots,N_i.
	\end{equation*}
	For $i=1,2,3$, the index value $j=0$ corresponds to the graph vertex $v_i$, and the index value $j=N_i$ corresponds to the vertex $v_4$ --- the center of the star. In this case the system (\ref{Wave_eqn_d33})--(\ref{Kirch_d_1mm33}) takes the form:
	\begin{align}
		u^i_{j,t+1}+u^i_{j,t-1}-u^i_{j+1,t}-u^i_{j-1,t}=0,\   j=1,\ldots,N_i-1, t=0,\ldots\label{Wave_eqn_c3}\\
		u^1_{N_1,t}=u^2_{N_2,t}=u^3_{N_3,t},\  t=1,\ldots \label{Cont_c3}\\
		-\frac{3}{2}u^1_{N_1,t+1}-\frac{3}{2}u^1_{N_1,t-1}  + u^1_{N_1-1,t}+u^2_{N_2-1,t}+u^3_{N_3-1,t}=0.\label{Kirch_c_1mm3}\\
		u^i_{0,t}=f^i_t,\ u^3_{0,t}=0, \ i=1,2,\ t=1,\ldots  \label{Bound_c3}\\
		u^i_{j,0}=u^i_{j,-1}=0,\ i=1,2,3,\ j=1,\ldots,N.\label{Init_c3}
	\end{align}
	Below we use the notation 
	$$
	g_t:= u^1_{N_1,t}=u^2_{N_2,t}=u^3_{N_3,t}=\left\{
	\begin{array}l 0,\quad \text{if $t<\min\{N_1,N_2,N_3\}$},\\
		\not=0\quad \text{if $t\geqslant\min\{N_1,N_2,N_3\}$}.
	\end{array}
	\right.
	$$
	Now it is easy to write down the solution of the forward problem to the 3-star graph using formulas obtained for the finite interval (\ref{solution_left}), (\ref{solution_right}). Namely,
	\begin{equation}
		u^1_{j,t}= u_{j,t}^{f^1-} + u_{j,t}^{g+},\quad  u^2_{j,t}= u_{j,t}^{f^2-} + u_{j,t}^{g+}, \quad u^3_{j,t}=  u_{j,t}^{g+}\label{solution_fp}.
	\end{equation}
	Since boundary controls are known, it remains to find $g_t$ to solve the forward problem on the discrete graph.
	Formula for the function $g_t$ can be derived using Table \ref{table3}, formula (\ref{solution_left}) and linearity of the system:
	\begin{equation}
		g_t = \frac{2}{3} (f^1_{t-N_1}+f^2_{t-N_2}).\label{solution_center}
	\end{equation}
	This equality is valid for $t<N_1+N_2+N_3+\min{N_1,N_2,N_3}-\max{N_1,N_2,N_3}$.
	
	\subsection{Shape control problem}
	The \emph{shape control} problem on a 3-star graph is to determine the boundary controls $f^1_t,f^2_t$ at $v_1$ and $v_2$ such that $u$, the solution of (\ref{Wave_eqn_c3})--(\ref{Init_c3}), satisfies the equalities  $u^i_{j,T}=\varphi^i_j$,
 $i=1,2,3$
and $j=1,...,N_i$,
 where 
	$\varphi=(\varphi^1,\varphi^2,\varphi^3)$ is a prescribed shape at time $t=T$. The solution of the shape control problem is based on the solution of the forward problem (\ref{solution_fp}). We will demonstrate that the
	optimal control time is
	$$
	T=\min\{L_1,L_2\},\quad L_1:=\max\{N_1+N_3,N_2\},\quad L_2:=\max\{N_1,N_2+N_3\}. 
	$$
	Since no control acts on $v_3$, edge $e_3$ should be controlled through $g_t$. Therefore,
	$$
 g_{T-j}=\left\{
	\begin{array}l   \varphi^3_{N_3-j},\quad j=0,\ldots,N_3,\\
		0,\quad j=N_3,\ldots,T.
  \end{array}
	\right.
  $$

	Note that $g_t$ can be controlled both from $v_1$ and $v_2$ (\ref{solution_center}).
	Let us consider the case when $L_1\leqslant L_2$.  Then we can find control
	$\tilde f^1$ on $v_1$ such that if $f^2=0$ on $v_2$ then $u^i_{N_i,t}=g_t$,
  $i=1,2,3$
and $t=0,....$:
	$$
	\tilde f^1_{t-N_1}=\frac{3}{2}g_t, \quad  t\leq N_1.
	$$
	The control we interested in has a form $f^1=\tilde f^1 + \hat f^1$ and $\hat f^1$ can be found as a solution to the equation
	\begin{equation}
		u^{\tilde f^1-}+u^{\hat f^1-}+u^{g+}=\varphi^1\quad \text{on $e_1$}.\label{eq_hat_f1}
	\end{equation}
	Such control $f^1$ and $f^2=0$ makes solution to (\ref{Wave_eqn_c3})--(\ref{Init_c3}) satisfying $u^i_{j,T}=\varphi^i_j$ $i=1,2$. It remains to satisfy the third relation $u^3_{j,T}=\varphi^3_j$ by appropriate choice of $f_2$. Exactly, $f^2$ should be a solution to the equation
	\begin{equation}
		u^{f^2-}+ u^{g+}=\varphi^2\quad \text{on $e_2$}.\label{eq__f2}
	\end{equation}
	Solving the equations (\ref{eq_hat_f1}) and (\ref{eq__f2}) with respect to $\hat f^1$ and $f^2$ we solve the shape control problem. The case $L_1>L_2$ can be considered similarly.
	
	For  solutions of equations  (\ref{eq_hat_f1}) and (\ref{eq__f2}) an explicit formulas can be written, we present
	ones in the case when $N_1=N_2=N_3$: 
	$$
	f^1_j=\left\{
	\begin{array}l   \frac{3}{2}\varphi^3_{j+1},\quad j=0,1,\ldots,N-1,\\
		\varphi^1_{2N-j-1} + \frac{1}{2}\varphi^3_{2N-j-1},\quad j=N,N+1,\ldots,2N-1,
		
	\end{array}
	\right.
	$$ 
	$$
	f^2_j=\left\{
	\begin{array}l   0,\quad j=0,1,\ldots,N-1,\\
		\varphi^2_{2N-j-1} -\varphi^3_{2N-j-1},\quad j=N,N+1,\ldots,2N-1.
	\end{array}
	\right.
	$$
	
	In a general situation  due to the different lengths of the edges of the graph, waves can reflect not only from the inner vertex, but also from the vertex $v_1$ or $v_2$. In such a case the formulas are more involved and contain terms taking into account these reflections. In our subsequent publications we plan to study control and inverse problems on more general discrete graphs and compare the results with their continuous analogs.

\section*{Acknowledgments}
	The research of Sergei Avdonin was  supported  in part by the National Science Foundation, grant DMS 1909869, and by Moscow Center for Fundamental and Applied Mathematics.
	The work of Victor Mikhaylov and Alexander Mikhaylov was supported by the RFBR grant 20-01-00627.
The research of Abdon Choque-Rivero was  supported
 by CONACYT Project A1-S-31524 and CIC-UMSNH, Mexico.




\subsection*{Conflict of interest}

The authors declare no potential conflict of interests.


\bibliography{wileyNJD-AMA}{}

\begin{thebibliography}{10}
\providecommand \doibase [0]{http://dx.doi.org/}%

\bibitem{LLS-book}
Lagnese J, Leugering G, Schmidt EJPG. {\it Modelling, Analysis, and Control of
  dynamical elastic multilink structures}.
\newblock Basel: Birkhauser .
\newblock 1994.
\newblock ISBN 9781461202738.

\bibitem{ref15}
Lagnese JE, Leugering G, Schmidt EJPG. Modelling of dynamic networks of thin
  thermoelastic beams. {\it Mathematical Methods in the Applied Sciences}
  1993\string; 16(5)\string: 327--358.

\bibitem{KIIK}
Kiik JC, Kurasov P, Usman M. On vertex conditions for elastic systems. {\it
  Physics Letters A} 2015\string; 379(34)\string: 1871-1876.

\bibitem{ref17}
Berkolaiko G, Ettehad M. Three-dimensional elastic beam frames: Rigid joint
  conditions in variational and differential formulation. {\it Studies in
  Applied Mathematics} 2022\string; 148(4)\string: 1586--1623.

\bibitem{GL-journal}
Gugat M, Leugering G. Global boundary controllability of the Saint-Venant
  system for sloped canals with friction. {\it Ann. Inst. H. Poincar\'e Anal.
  Non Lineaire} 2009\string; 26(1)\string: 257--270.

\bibitem{ref1}
Steinbach MC. On pde solution in transient optimization of gas networks. {\it
  Journal of Computational and Applied Mathematics} 2007\string; 203(2)\string:
  345--361.

\bibitem{ref2}
Domschke P, Kolb O, Lang J. Adjoint-based error control for the simulation and
  optimization of gas and water supply networks. {\it Applied Mathematics and
  Computation} 2015\string; 259\string: 1003--1018.

\bibitem{ref3}
Egger H, Schobel-Krohn L. Chemotaxis on networks: analysis and numerical
  approximation. {\it ESAIM: Mathematical Modelling and Numerical Analysis}
  2020\string; 54(4)\string: 1339--1372.

\bibitem{ref4}
Bretti G, Natalini R, Ribot M. A hyperbolic model of chemotaxis on a network: a
  numerical study. {\it ESAIM: Mathematical Modelling and Numerical Analysis -
  Modélisation Mathématique et Analyse Numérique} 2014\string; 48(1)\string:
  231--258.

\bibitem{ref5}
Borsche R, Klar A, Pham TNH. Nonlinear flux-limited models for chemotaxis on
  networks. {\it Networks and Heterogeneous Media}
  2017(1556-180L2017\_3\_381)\string: 381--401.

\bibitem{ref6}
Herty M, Mohring J, Sachers V. A new model for gas flow in pipe networks. {\it
  Mathematical Methods in the Applied Sciences} 2011\string; 33(7)\string:
  845--855.

\bibitem{BCd-book}
Bastin G, Coron JM, Novel dB. Using hyperbolic systems of balance laws for
  modeling, control and stability analysis of physical networks. In: Lecture
  notes for the Pre-Congress Workshop on Complex Embedded and Networked Control
  Systems 17th IFAC World Congress. IFAC. ; 2008; Seoul, Korea.

\bibitem{Hante-book}
Hante FM, Leugering G, Martin A, Schewe L, Schmidt M. {\it Challenges in
  optimal control problems for gas and fluid flow in networks of pipes and
  canals: From modeling to industrial applications}.
\newblock Singapore: 122, Ind. Appl. Math., Springer .
\newblock 2017.
\newblock ISBN 9048154464.

\bibitem{ABG-journal}
Al\`i G, Bartel A, G{\"u}nther M. Parabolic Differential-Algebraic Models in
  Electrical Network Design. {\it Multiscale Model. Simul.} 2005\string;
  4(3)\string: 813--838.

\bibitem{ref14}
Cheng X, Scherpen JMA. Clustering approach to model order reduction of power
  networks with distributed controllers. {\it Advances in Computational
  Mathematics volume} 2018\string; 44\string: 1917--1939.

\bibitem{CGHS-journal}
Colombo RM, Guerra G, Herty M, Schleper V. Optimal control in networks of pipes
  and canals. {\it SIAM J. Control Optim.} 2009\string; 48(3)\string:
  2032--2050.

\bibitem{ref7}
Garavello M, Piccoli B. {\it Traffic Flow on Networks, volume 1 of AIMS Series
  on Applied Mathematics}.
\newblock Springfield, MO: American Institute of Mathematical Sciences .
\newblock 2006.

\bibitem{ref9}
Oppenheimer SF. A convection-diffusion problem in a network. {\it Applied
  Mathematics and Computation} 2000\string; 112(2)\string: 223--240.

\bibitem{ref10}
Herty M, Ringhofer C. Averaged kinetic models for flows on unstructured
  networks. {\it Kinetic and Related Models} 2011\string; 4(4)\string:
  1081--1096.

\bibitem{ref11}
García L, Barreiro-Gomez J, Escobar E, Téllez D, Quijano N, Ocampo-Martinez
  N. Modeling and real-time control of urban drainage systems: A review. {\it
  Advances in Water Resources} 2015\string; 85\string: 120--132.

\bibitem{H-book}
Hurt NE. {\it Mathematical physics of quantum wires and devices. From spectral
  resonances to Anderson localization}.
\newblock Dordrecht: Kluwer Academic .
\newblock 2000.
\newblock ISBN 9048154464.

\bibitem{JR-book}
Joachim C, Roth S. {\it Atomic and Molecular Wires}.
\newblock Dordrecht: Springer .
\newblock 1997.
\newblock ISBN 9780792346289.

\bibitem{MP-journal}
Melnikov YuB, Pavlov BS. Two-body scattering on a graph and application to
  simple nanoelectronic devices. {\it J. Math. Phys.} 1995\string;
  36(6)\string: 2813--2825.

\bibitem{ref12}
Duca A. Bilinear quantum systems on compact graphs: Well-posedness and global
  exact controllability. {\it Automatica} 2021\string; 123\string: 123:109324.

\bibitem{KS1-journal}
Kottos T, Smilansky U. Quantum chaos on graphs. {\it Phys. Rev. Lett.}
  1997\string; 79(24)\string: 4794--4797.

\bibitem{Hirt1974}
Adam S, Hwang E, Galits VM, Das~Sarma S. A self-consistent theory for graphene
  transport. {\it Proceeding of the National Academy of Sciences} 2007\string;
  104(47)\string: 18392--18397.

\bibitem{ref8}
Kuchment P, Post O. On the spectra of carbon nano-structures. {\it Commun.
  Math. Phys} 2007\string; 275\string: 805--826.

\bibitem{BC-journal}
Bell J, Cracium G. A distributed parameter identification problem in neuronal
  cable theory models. {\it Math. Biosci.} 2005\string; 194(1)\string: 1--19.

\bibitem{AB-journal}
Avdonin S, Bell J. Determining a distributed conductance parameter for a
  neuronal cable model defined on a tree graph. {\it Inverse Problems and
  Imaging} 2016\string; 9(3)\string: 645--659.

\bibitem{ref13}
Du B, Lian X, Cheng X. Partial differential equation modeling with Dirichlet
  boundary conditions on social networks. {\it Boundary Value Problems}
  2018\string; 50\string: 2035--2052.

\bibitem{ref19}
Solomon J. PDE approaches to graph analysis. {\it ArXiv} 2015\string;
  abs/1505.00185.

\bibitem{ref20}
Mercier D, Régnier V. Spectrum of a network of Euler-Bernoulli beams. {\it
  Journal of Mathematical Analysis and Applications} 2008\string;
  337(1)\string: 174--196.

\bibitem{ref21}
Kottos T, Smilansky U. Periodic orbit theory and spectral statistics for
  quantum graphs. {\it Annals of Physics} 1999\string; 274(1)\string: 76--124.

\bibitem{ref22}
{von Below} J. A characteristic equation associated to an eigenvalue problem on
  c2-networks. {\it Linear Algebra and its Applications} 1985\string;
  71\string: 309--325.

\bibitem{ref23}
Band R, Lévy G. Quantum graphs which optimize the spectral gap. {\it Annales
  Henri Poincaré} 2017\string; 18\string: 3269--3323.

\bibitem{ref24}
Laurent M, Piovesan T. Conic approach to quantum graph parameters using linear
  optimization over the completely positive semidefinite cone. {\it SIAM
  Journal on Optimization} 2015\string; 25(4)\string: 2461--2493.

\bibitem{BK-book}
Berkolaiko G, Kuchment P. {\it Introduction to Quantum Graphs (Mathematical
  Surveys and Monographs vol. 186)}.
\newblock Providence, RI: American Mathematical Society .
\newblock 2013.
\newblock ISBN 0821892118.

\bibitem{ref25}
Domschke P, Dua A, Stolwijk J, Lang J, Mehrmann V. Adaptive refinement
  strategies for the simulation of gas flow in networks using a model
  hierarchy. {\it Electronic Transactions on Numerical Analysis} 2018\string;
  48\string: 97--113.

\bibitem{ref26}
Egger H, Philippi N. A hybrid discontinuous galerkin method for transport
  equations on networks. In: In Robert Klofkorn, Eirik Keilegavlen, Florin A.
  Radu, and J{\"u}rgen Fuhrmann, editors, Finite Volumes for Complex
  Applications IX - Methods, Theoretical Aspects, Examples. Dynamic Publishers.
  ; 2020; Springer International Publishing\string: 487--495.

\bibitem{ref27}
Grundel S, Herty M. Hyperbolic discretization of simplified Euler equation via
  Riemann invariants. {\it Applied Mathematical Modelling} 2022\string;
  106\string: 60--72.

\bibitem{ref28}
Grundel S, Hornung N, Roggendorf S. Numerical aspects of model order reduction
  for gas transportation networks. In: In Slawomir Koziel, Leifur Leifsson, and
  Xin-She Yang, editors, Simulation Driven Modeling and Optimization. Dynamic
  Publishers. ; 2016; Springer International Publishing\string: 1--28.

\bibitem{ref29}
Gugat M, Herty M, Schleper V. Flow control in gas networks: Exact
  controllability to a given demand. {\it Mathematical Methods in the Applied
  Sciences} 2011\string; 34(7)\string: 745--757.

\bibitem{ref30}
Pesenson I. Polynomial splines and eigenvalue approximations on quantum graphs.
  {\it Journal of Approximation Theory} 2005\string; 132(2)\string: 203--220.

\bibitem{ref31}
Wybo WA, Boccalini D, Torben-Nielsen B, Gewaltig MO. Asparse reformulation of
  the green's function formalism allows efficient simulations of morphological
  neuron models. {\it Neural Computation} 2015\string; 27(12)\string:
  2587--2622.

\bibitem{ref32}
Arioli M, Benzi M. A finite element method for quantum graphs. {\it IMA Journal
  of Numerical Analysis} 2017\string; 38(3)\string: 1119--1163.

\bibitem{ref33}
García L, Barreiro-Gomez J, Escobar E, Téllez D, Quijano N, Ocampo-Martinez
  C. Nonoverlapping domain decomposition for optimal control problems governed
  by semi-linear models for gas flow in networks. {\it Control and Cybernetics}
  2017\string; 46(3)\string: 191--225.

\bibitem{ref34}
Stoll M, Winkler M. Optimal Dirichlet control of partial differential equations
  on networks. {\it Electronic transactions on Numerical Analysis} 2021\string;
  54\string: 392--419.

\end{thebibliography}
\bibliographystyle{plain}

\end{document}